\documentclass[reqno,12pt]{amsart}
\usepackage[english]{babel}
\usepackage{amsmath, amssymb, amsfonts}
\usepackage{epsfig} 
\usepackage{hyperref}

\theoremstyle{definition}
\newtheorem{defin}{Definition}

\theoremstyle{plain}
\newtheorem{lemma}{Lemma}
\newtheorem{theorem}{Theorem}
\newtheorem{sled}{Corollary}
\newcommand{\NN}{{\mathbb N}}
\newcommand{\RR}{{\mathbb R}}
\newcommand{\ZZ}{{\mathbb Z}}

\begin{document}
\thispagestyle{empty}

\title{On the properties of functions\\of the Takagi power class}

\author{O.E.~Galkin, S.Yu.~Galkina, A.A.~Tronov}

\address{Galkin Oleg Evgenievich,
\newline\hphantom{iii} National Research University ``Higher School of Economics'', Russian Federation}
\email{olegegalkin@ya.ru}

\address{Galkina Svetlana Yurievna,
\newline\hphantom{iii} National Research University ``Higher School of Economics'', Russian Federation}
\email{svetlana.u.galkina@mail.ru}

\address{Tronov Anton Alexandrovich,
\newline\hphantom{iii} National Research University ``Higher School of Economics'', Russian Federation}
\email{tronovaa@yandex.ru}

\maketitle

\tableofcontents
 
{\small
\begin{quote}
\medskip

\noindent{\bf Abstract.} 
By construction, functions of Takagi power class are similar to Takagi's continuous
nowhere differentiable function. These functions have one real parameter $p>0$. 
They are defined on the real axis by the series $S_p(x) = \sum_{n=0}^\infty (T_0(2^nx)/2^n)^p$,
where $T_0(x)$ is the distance from $x\in\RR$ to the nearest integer.
We study such properties of functions $S_p(x)$ as continuity, generalized H\"older condition, global maxima and differentiability at the points $x=1/3$, $x=2/3$, for various  values of $p$.
\medskip

\end{quote} 
}

\section{Introduction}
\label{sectVved}
This article is devoted to the study of functions from the Takagi power class.
These functions are similar in construction to the continuous but nowhere differentiable Takagi function described in 1903 (\cite{Takagi1903}). They have one real parameter $p>0$ and are defined as follows.


\begin{defin}
We call \emph{the Takagi power class} a family of functions $S_p(x)$ having one real parameter $p>0$ and given on the real axis $\RR$ by the equality
\begin{equation}\label{eqDefSp}
S_p(x) = \sum_{n=0}^{\infty}\left( \frac{T_{0}(2^{n} x)}{2^{n}} \right)^p, \quad x\in\RR,
\end{equation}
where $T_0(x) = |x-[x+1/2]| = |\{x+1/2\}-1/2| = \rho(x, {\mathbb Z})$
is the distance between the point $x$ and the nearest integer point, 
$[y]$ is the integer part (respectively $\{y\}$ is the fractional part) of the number $y\in\RR$.

The function $T_{0}$ can also be defined as follows:
\begin{equation} \label{eqDefT0}
T_{0}(x) =  
\begin{cases}
x-n     & \text{for } x \in [n, n + 1/2], n\in\ZZ;\\
1-(x-n) & \text{for } x \in [n + 1/2, n + 1], n \in \ZZ. 
\end{cases}
\end{equation}
\end{defin}
Note that for $p=1$, the function $S_p(x)$ coincides with the Takagi function.
Hata and Yamaguti~\cite[Sec.~2]{HataYamaguti84} replaced the sequence of coefficients $\{1/2^n\}$ 
in the definition of the Takagi function with an arbitrary sequence of constants
and got a new family of functions, calling it \emph{Takagi class}.

The real functions $T_v(x)$ which have on the real parameter $v\in(-1;1)$ and can be defined by the formula
\begin{equation} 
T_v(x) = \sum\limits_{n=0}^\infty v^n T_0(2^nx), \quad x\in {\mathbb R}.
\end{equation}
form a narrower class, namely {\it exponential Takagi-Landsberg class} 
(see \cite{Landsberg1908}, \cite{MishSchied2019}, \cite{Han2019}, \cite{HanSchied2020}, \cite{GalkinArx2020}).
Note that for $v=0$, the function $T_v(x)$ matches $T_0(x)$, and for $v=1/2$ it matches with the Takagi function $T(x)$.

We can say that the problem of finding extremes for $T_v$ was set by J.~Tabor and J.~Tabor (see~\cite[Problem 1.2, p.~731]{Tabor2}.
In order to accurately estimate continuous semi-convex functions, they introduced the functions $\omega_p$, which can be set as $\omega_p(x) = 2\cdot T_{1/2^p}(x)$, and obtained the formula for global maxima of the functions $\omega_{p_n}$ on $[0;1]$ for one specific sequence $\{p_n\}$ (see~\cite[Theorem~3.1]{Tabor2}).

Here are some more results about the extremes of the $T_v$ functions.

1) In the case of $v=1/2$ Kahane (see~\cite[Lieu~1]{Kahane}) found points of local and global extremes for $T_v (x)=T(x)$.

2) Further results on local extremes and sets of the Takagi function level can be found in the reviews~\cite{AllaartKawamuraSurv} and~\cite{Lagarias}.

3) In~\cite[Theorem 4]{Galkin2015} it is proved that for $v\in[-1/2;1/4]$, the point $1/2$ is the point of the global maximum of the function $T_v$ on $[0;1]$, and for $v\in(-1;-1 / 2)\cup(1/4;1)$ is not.

4) In the case of $v=-1/2$, it follows from Allaart's remark~\cite[Remark 5.6, p.28]{Allaart} 
that the set of minimum points of the function $T_v$ is a Cantor type set obtained by removing the ``middle half'', and therefore is uncountable.

5) Han and Schied in \cite{Han2019},~\cite{HanSchied2020} gave a new approach to characterizing and computing 
the set of global maximizers and minimizers of the functions in the Takagi class and, in particular, 
in the exponential Takagi–Landsberg class. They showed that the function $T_v$ has a unique maximizer in $[0,1/2]$ 
if and only if there does not exist a Littlewood polynomial that has $2v$ as a certain type of root, called step root. 
Their general results lead to explicit and closed-form expressions for the maxima 
of the exponential Takagi–Landsberg functions with $v\in(-1,1/4]\cup(1/,12)$.

Similar results for maximums and minimums of functions in exponential Takagi–Land\-s\-berg class were obtained. in~\cite{GalkinArx2020}.

The Takagi function and its generalizations are used in various fields of mathematics: mathematical analysis, probability theory, number theory, and others. A large number of publications are dedicated to these functions, and this number continues to increase. Lots of interesting results and links 
available in the reviews~\cite{AllaartKawamuraSurv} and~\cite{Lagarias}.
See also more recent papers, for example \cite{Schied2015}, \cite{JaerHir2017}, \cite{Mishura2018}, \cite{Ferrera2019}, \cite{GalkinDelange2020}, \cite{MonroeArx2020}, \cite{Yu2020},~\cite{MonroeArx2021}.

In this article, for different values of the parameter $p$, we study such properties of functions $S_p(x)$ as
continuity, generalized H\"older condition, global maxima and differentiability at points $x=1/3$, $x=2/3$.

Let us briefly describe the structure and main results of the work.
The work consists of six sections.

{\it Section~\ref{sectVved}} is an introduction.

{\it Section~\ref{sectGrafsSp}} contains the graphs of functions $y=S_p(x)$ in Takagi power class for $p=0.2$, $p=0.4$, $p=0.7$ and $p=1$.

{\it In~section~\ref{sectNeprSp}} we prove that the functions $S_p(x)$ is everywhere defined in $\RR$, continuous and has a period of $1$.

{\it In~section~\ref{sectHolderSp}} we prove that for every $p\in(0;1)$ the functions $S_p(x)$ 
satisfy the generalized H\"older condition with the exponent $p$.

{\it In~section~\ref{sectFunkUrSp}} we get some functional equation that the functions $S_p(x)$ satisfy for every $p\in(0;1)$.

{\it The~section~\ref{sectMaxSp}} is devoted to finding the global maxima of the functions $S_p(x)$ and the sets where they are reached. Moreover, here we prove that for any $p\in(0;1)$ the function $S_p(x)$ has no derivative at points $x=1/3$ and $x=2/3$.

\section{Graphs of functions $y=S_p(x)$ of the Takagi power class}
\label{sectGrafsSp}

In this section, we present sketches of graphs of functions $y=S_p(x)$
from the Takagi power class for $p=0.2$, $p=0.4$, $p=0.7$ and $p=1$.
Two vertical lines indicate two points of the global maximum, $x=1/3$ and $x=2/3$.

\newlength{\VisotaRis}
\VisotaRis = 0.42\textheight
\begin{figure}[h] 
\centering
\includegraphics[height=\VisotaRis]{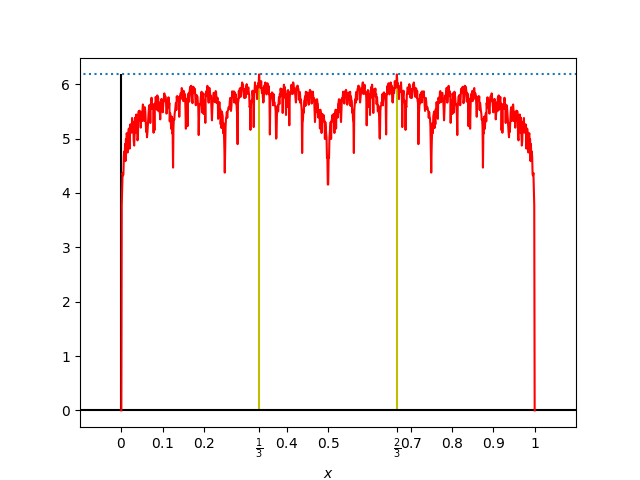}
\caption{Graph of the function $y=S_p(x)$ for $p = 0.2$.}
\label{a_02}
\end{figure}

\begin{figure}[h] 
\centering
\includegraphics[height=\VisotaRis]{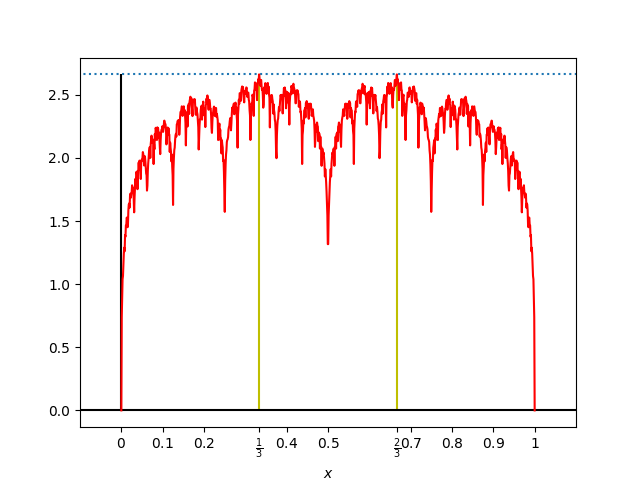}
\caption{Graph of the function $y=S_p(x)$ for $p = 0.4$.}
\label{a_04}
\end{figure}

\begin{figure}[h] 
\centering
\includegraphics[height=\VisotaRis]{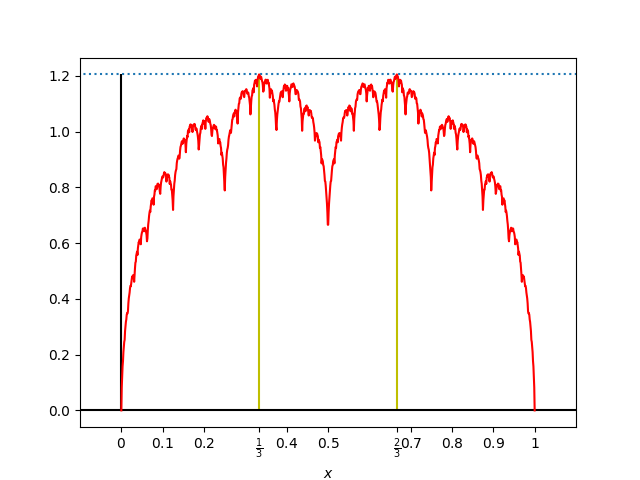}
\caption{Graph of the function $y=S_p(x)$ for $p = 0.7$.}
\label{a_07}
\end{figure}

\begin{figure}[h] 
\centering
\includegraphics[height=\VisotaRis]{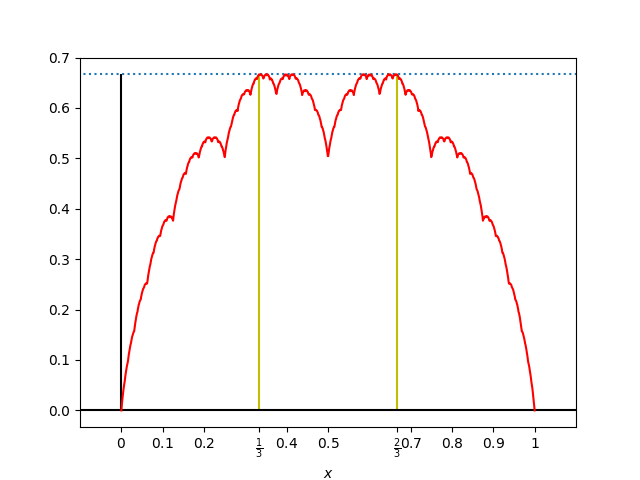}
\caption{Graph of the function $y=S_p(x)$ for $p = 1$.}
\label{a_1}
\end{figure}    

\section{Continuity of functions of the Takagi power class}
\label{sectNeprSp}
In this section, we present  proof of Takagi power class  functions continuity.

\begin{theorem}
For any $p>0$, the function $S_p(x)$ from the Takagi power class is everywhere defined on $\RR$, is continuous, has a period of $1$ and has the following symmetry property:
\begin{equation}\label{eqSimmSp}
S(x) = S(q-x) \quad\text{for all }q\in\ZZ\text{ and all}x\in\RR.
\end{equation}
In addition, for all $p>0$ and $x\in\RR$, the inequality $|S_p(x)|\leq 1/(2^p-1)$ holds.
\end{theorem}
\begin{proof}
1) First we will use the Weierstrass M-test to prove the uniform convergence of the functional series~\eqref{eqDefSp}.

Obviously, for all $p>0$ and $x\in\RR$, the inequality 
$|T_0(x)|\leq 1/2$ holds.
Therefore
$$
\left|\left( \frac{T_{0}(2^{n} x)}{2^{n}} \right)\right|^p \leq\left(\frac{1/2}{2^{n}} \right)^p =
\frac{1}{2^{(n+1)p}}.
$$
The series $\sum_{n=0}^{\infty}1/2^{(n+1)p}$ converges, and its sum is equal to $1/(2^p-1)$.
Hence, the series~\eqref{eqDefSp} defining $S_p(x)$ converges uniformly on all $\mathbb{R}$, 
with $S_p(x)\leq 1/(2^p-1)$ for all $x\in\RR$.

2) From the uniform convergence of the series~\eqref{eqDefSp} and the continuity of its terms,
according to the theorem on the continuity of the sum of a functional series, the continuity of the function
$S_p(x) = \sum_{n=0}^{\infty}\big(T_{0}(2^{n}x)/2^{n}\big)^p$
for all $x\in\RR$ follows.

3) The periodicity of the function $S_p(x)$ follows from the formula~\eqref{eqDefSp} and the presence of the period $1$ of the function $T_0$.
\end{proof}

\section{Generalized H\"older condition for functions of the Takagi power class}
\label{sectHolderSp}
In this section, we present  proof of H\"older condition for functions of the Takagi power class, based on some proved lemmas.
\begin{lemma} \label{lemHolder1}
For any real numbers $x,y\in\mathbb{R}$ and any $p>0$, the inequality holds
\begin{equation} \label{eqT0xT0yOgr}
|T_{0}^p(x) - T_{0}^p(y)| \leq\frac{1}{2^p}.
\end{equation}
\end{lemma}
\begin{proof}
Since $T_{0}(x) \in [0, 1/2]$ for all $x\in\RR$, then $T_{0}^p(x)\in[0,1/2^p]$. 
Therefore, for all $x,y\in\RR$, the inequality $|T_{0}^p(x) - T_{0}^p(y)| \leq 1/2^p$ holds.
\end{proof}

\begin{lemma} \label{lemHolder2}
For any real numbers $x, y\in\mathbb{R}$ and any $p\in (0;1]$, the inequality holds
\begin{equation} \label{eqHoldT0p}
|T_{0}^p(x) - T_{0}^p(y)| \leq|x - y|^p
\end{equation}
\end{lemma}
\begin{proof}
1) Let's first prove that the following inequality is true:
\begin{equation} \label{eqABP}
|a^p-b^p| \leq |a-b|^p \quad\text{for all } a\geq0, b\geq0\text{ and all } p\in(0;1].
\end{equation}
Assume that $a\geq b$.
The case $a = 0$ is trivial.

Let $a > 0$. Dividing both parts of the inequality~\eqref{eqABP} by $a^p>0$, we get an equivalent inequality: 
$|1 -(b/a)^p| \leq|1 - b/a|^p$.
Since $b/a\in[0;1]$ and $p\in(0;1]$, this inequality follows from the double inequality
$1 -(b/a)^p \leq 1 - b/a \leq (1 - b/a)^p$.

It is easy to see that in the case $a<b$, the inequality~\eqref{eqABP} also holds, due to symmetry: 
$|a^p - b^p| = |b^p - a^p| \leq |b - a|^p = |a - b|^p$.

2) For all $x,y\in\RR$, let us prove the inequality 
\begin{equation} \label{eqT0xy}
|T_{0}(x) - T_{0}(y)|\leq|x-y|.
\end{equation}
Since the function $T_0(x)$ is continuous on $\RR$ and there is a derivative $T_0'(s)=\pm1$ at all points $s\in\RR$, except for points of the form $s=q/2$, where $q\in\ZZ$, the relations are true
$$
|T_{0}(x) - T_{0}(y)| = \left|\int_{x}^{y}T_0'(s)ds\right| \leq 
\left|\int_{x}^{y}|T_0'(s)|ds\right| = |x-y|.
$$

3) Based on the results of points 1) and 2) of our proof, we prove the inequality~\eqref{eqHoldT0p} for all $x,y\in\RR$.

Let's put $a=T_0(x)$ and $b=T_0(y)$. Then $a\geq0$ and $b\geq0$, so the inequality~\eqref{eqABP} implies an estimate of $|T_{0}^p(x) - T_{0}^p(y)| \leq|T_{0}(x) - T_{0}(y)|^p$.
Hence, raising inequality~\ref{eqT0xy} to the power $p$, we deduce the desired relation~\eqref{eqHoldT0p}:
$$
|T_{0}^p(x) - T_{0}^p(y)| \leq |T_{0}(x) - T_{0}(y)|^p \leq |x-y|^p.
$$
\end{proof}

\begin{lemma} \label{lemTehnOc}
For any $p\in[0;1]$, the following inequality holds:
$p\cdot\ln{2}\cdot 2^p \geq 2^p - 1$.
\end{lemma}
\begin{proof}
The inequality being proved is equivalent to the inequality
$2^p(1-p\cdot\ln{2}) \leq 1$.

Let's define the function $g(p) = 2^p(1-p\cdot\ln{2})$.
This function does not increase on the segment $[0;1]$, since
$$
g'(p) = 2^p( (1-p\cdot\ln{2})\ln{2} - \ln{2} ) = -p\cdot2^p\ln^p{2} \leq 0.
$$
Also, $g(0)=$1. Therefore, $g(p) = 2^p(1-p\cdot\ln{2}) \leq 1$ for all $p\in[0;1]$.
\end{proof}

\begin{defin}
We will say that the function $f\colon\RR\to\mathbb{R}$ satisfies the H\"older condition with the exponent $p$ and with the constant $C$, if for all $x,y\in\RR$ the inequality holds
$$
|f(x)-f(y)| \leq C  |x-y|^p. 
$$
\end{defin}

\begin{defin}
We will say that the function $f\colon\RR\to\mathbb{R}$ satisfies the generalized H\"older condition with the exponent $p$ and with the constant $C$, if for all different $x,y\in\RR$ such that $|x-y|\leq1/2$, the inequality holds
$$
|f(x)-f(y)| \leq C  |x-y|^p\!\cdot\log_2\frac{1}{|x-y|}. 
$$
\end{defin}

\begin{theorem}
For any $p\in(0;1]$ function 
$S_p(x) =  \sum_{k=0}^{\infty}(T_{0}(2^{k}x)/2^{k})^p$
of Takagi power class satisfies the generalized H\"older condition with the exponent $p$ and constant
$C = 1/p \cdot \log_{2}\big(p\cdot\ln{2}/(2^p-1)\big) + 2^p/(p\cdot\ln{2}) +1$.
\end{theorem}    
\begin{proof}
Suppose there are two distinct points $x,y\in\RR$ that holds $|x-y|\leq1/2$. Let's fix them.
Then the following relations hold:
\begin{equation} \label{eqSpxSpy}
\begin{split}
|S_p(x) - S_p(y)| = 
\left| \sum_{k=0}^{\infty}\left(\frac{T_{0}(2^k x)}{2^{k}} \right)^p - 
\sum_{k=0}^{\infty}\left(\frac{T_{0}(2^k y)}{2^{k}} \right)^p \right| = 
\left| \sum_{k=0}^{\infty}\frac{T_{0}^p(2^k x) - T_{0}^p(2^k y)}{2^{kp}} \right| \leq\\
\leq \sum_{k=0}^{\infty}\frac{ | T_{0}^p(2^k x) - T_{0}^p(2^k y) |}{2^{kp}} =
\sum_{k=0}^{n-1}\frac{| T_{0}^p(2^k x) - T_{0}^p(2^k y)|}{2^{kp}} +
\sum_{k=n}^{\infty}\frac{ | T_{0}^p(2^k x) - T_{0}^p(2^k y) |}{2^{kp}}.
\end{split}
\end{equation}
Let's estimate the terms of the first sum in the right side~\eqref{eqSpxSpy} using inequality~\eqref{eqHoldT0p} of lemma~\ref{lemHolder2}: 
$$
|T_{0}^p(2^k x) - T_{0}^p(2^k y)| \leq |2^k x - 2^k y|^p.
$$ 
Let's estimate the terms of the second sum in the right side\eqref{eqSpxSpy} using inequality~\eqref{eqT0xT0yOgr} of lemma~\ref{lemHolder1}: 
$$
|T_{0}^p(2^k x) - T_{0}^p(2^k y)| \leq \frac1{2^p}.
$$
Using these two estimates, we are able to get from~\eqref{eqSpxSpy} statement below:
\begin{equation} \label{eqMainIneq}
|S_p(x) - S_p(y)| \leq
\sum_{k=0}^{n-1}\frac{|2^k x - 2^k y|^p}{2^{kp}} + \sum_{k=n}^{\infty}\frac{ 1/2^ p }{2^{kp}} = 
n \cdot | x -  y |^p + \frac{1}{2^{ p n}  (2^p - 1)}.
\end{equation}
Let's select $n \in \mathbb{N}$, so that to decrease right side.
For this goal we should find function minimum $f(t) = t\cdot|x-y|^p + 1/2^{pt}(2^p - 1)$.

It's derivative is:
$$
f'(t) = |x-y|^p - \frac{p \cdot \ln{2}}{2^{pt}(2^p - 1)}.
$$
From the condition $f'(t)=0$ we can find the point of maximum:
\begin{equation} \label{eqt0}
t_0 = \frac{1}p \cdot \log_{2}\Big( \frac{ p \cdot \ln{2}}{(2^p - 1) \cdot |x-y|^p} \Big).
\end{equation}
As $|x-y|\leq 1/2$ then $t_0 \geq 1/p \cdot \log_2\big(p\cdot\ln2\cdot 2^p/(2^p - 1)\big)$.
Therefore $t_0\geq 0$ according to the lemma~\ref{lemTehnOc}.

Let $n_0 = [t_0] = \big[1/p \cdot \log_2\big(p\cdot\ln2/((2^p - 1)\cdot|x-y|^p)\big)\big]$.
Then $n_0\geq 0$ and  inequality holds $t_0-1 < n_0 \leq t_0$.
According to the statement and~\eqref{eqt0} it follows that
\begin{equation} \label{eqLeN0}
\frac{1}p \cdot \log_{2}\Big( \frac{ p \cdot \ln{2}}{(2^p - 1) \cdot |x-y|^p} \Big) - 1 < n_0
\end{equation}
and
\begin{equation} \label{eqN0Le}
n_0 \leq \frac{1}p \cdot \log_{2}\Big( \frac{ p \cdot \ln{2}}{(2^p - 1) \cdot |x-y|^p} \Big).
\end{equation}
According to~\eqref{eqMainIneq}:
\begin{equation} \label{eqSpxSpyN0}
|S_p(x) - S_p(y)|  \leq n_0\cdot|x-y|^p + \frac{1}{2^{pn}(2^p - 1)}.
\end{equation}
According to~\eqref{eqLeN0} follows that
\begin{equation} \label{eqLeN0sled}
\frac{1}{2^{pn_0}(2^p - 1)} < \frac{2^p|x-y|^p}{p\cdot\ln{2}}.
\end{equation}
Multiplying both parts of inequality~\eqref{eqN0Le} on $|x-y|^p$, we have got:
\begin{equation} \label{eqN0LeSled}
n_0 \cdot | x -  y |^p \leq\left( \frac{1}p \cdot \log_{2}\frac{p\cdot\ln{2}}{2^p - 1} + 
\log_{2}\frac{1}{|x-y|} \right) \cdot |x-y|^p.
\end{equation}
According~\eqref{eqLeN0sled} and~\eqref{eqN0LeSled} and formula~\eqref{eqSpxSpyN0}, we get the estimate 
$$
|S_p(x) - S_p(y)| \leq
\left( \frac{1}p \cdot \log_{2}\frac{ p\cdot\ln{2}}{2^p-1} + \frac{2^p}{p\cdot\ln{2}} + \log_{2}\frac{1}{|x-y|} \right) \cdot |x-y|^p.
$$
Let's denote $b = 1/p \cdot \log_{2}\big(p\cdot\ln{2}/(2^p-1)\big) + 2^p/(p\cdot\ln{2})$. 
This case we can rewrite previous estimate as
\begin{equation} \label{eqSpxSpyB}
|S_p(x) - S_p(y)| \leq \left( b + \log_{2}\frac{1}{|x-y|} \right) \cdot |x-y|^p.
\end{equation}
As $|x-y|\leq 1/2$, then $\log_2(1/|x-y|) \geq 1$.
Hence from~\eqref{eqSpxSpyB} the inequality below holds
$$
|S_p(x) - S_p(y)| \leq (b+1)\cdot|x-y|^p\log_{2}\frac{1}{|x-y|} =
C\cdot|x-y|^p\log_{2}\frac{1}{|x-y|},
$$
where 
$$
C = b + 1 = \frac{1}p \cdot \log_{2}\left(\frac{p\cdot\ln{2}}{2^p -1}\right) + \frac{2^p}{p\cdot\ln{2}} +1.
$$
\end{proof}

\section{Functional equation for functions from the Takagi power class}
\label{sectFunkUrSp}
The section contains proof and important particular cases of functional equation for Takagi power class  functions.
\begin{theorem}
For any $p>0$ and $m\in\NN$ functions $S_p(x)$  of Takagi power class satisfy the functional equation
\begin{equation} \label{eqGenFuncUr}
S_p(x) = \sum_{k=0}^{m-1}\left(\frac{T_{0}(2^k x)}{2^{k}} \right)^p + \frac{S_p(2^{m} x)}{2^{mp}}, \quad x\in\RR.
\end{equation} 
\end{theorem} 
\begin{proof}

According to definition~\eqref{eqDefSp} of $S_p(x)$:
$$
S_p(2^{m}x) = \sum_{n=0}^{\infty}\left(\frac{T_{0}(2^n\cdot 2^{m}x)}{2^n} \right)^p = \sum_{n=0}^{\infty}\left(\frac{T_{0}(2^{n+m}x)}{2^n} \right)^p.
$$
Dividing both parts by $2^{mp} > 0$, we get:
$$
\frac{S_p(2^{m}x)}{2^{mp}} = \sum_{n=0}^{\infty}\left(\frac{T_{0}(2^{n+m}x)}{2^{n+m}} \right)^p = 
\sum_{k=m}^{\infty}\left(\frac{T_{0}(2^{k}x)}{2^{k}} \right)^p = 
S_p(x) - \sum_{k=0}^{m-1}\left(\frac{T_{0}(2^k x)}{2^{k}} \right)^p.
$$ 
From this required equation~\eqref{eqGenFuncUr} follows.
\end{proof}

\begin{sled}
For any $p>0$ functions $S_p(x)$ of Takagi power class satisfy the functional equations
\begin{equation} \label{eqFuncUrM1}
S_p(2x) = 2^p\big(S_p(x) - T_{0}^p(x)\big), \quad x\in\RR       
\end{equation} 
and
\begin{equation} \label{eqFuncUrM2}
S_p(4x) = 4^p\big(S_p(x) - T_{0}^p(x)\big)- 2^pT_{0}^p(2x), \quad x\in\RR.      
\end{equation} 
\end{sled}
\begin{proof}
The functional equation~\eqref{eqFuncUrM1} follows from~\eqref{eqGenFuncUr} at $m=1$.
The functional equation~\eqref{eqFuncUrM2} follows from~\eqref{eqGenFuncUr} at $m=2$.
\end{proof}

\begin{lemma} \label{lemSp121315}
For any $p>0$ the values of the function $S_p(x)$ at the points $x=1/2$, $x=1/3$, $x=1/5$ and $x=2/5$ can be calculated by the following formulas:
$S_p(1/2)=1/2^p$, $S_p(1/3)=2^p/(6^p-3^p)$, $S_p(1/5)=2\cdot4^p/(20^p-5^p)$ and $S_p(2/5)=(8^p-2^p)/(20^p-5^p)$.
\end{lemma}
\begin{proof}
$ $

1) The equality $S_p(1/2)=1/2^p$ follows from the formulas \eqref{eqDefSp} and~\eqref{eqDefT0}.

2) Substituting the value $x=1/3$ into the formula~\eqref{eqFuncUrM1}, we get the following equation
$S_p(2/3) = 2^p\big(S_p(1/3) - T_0^p(1/3)\big)$.
Hence, since $T_0^p(1/3)=1/3$ by virtue of~\eqref{eqDefT0} and $S_p(2/3)= S_p(1/3)$ by virtue of the symmetry property~\eqref{eqSimmSp}, the following equality is true: $S_p(1/3) = 2^p\big(S_p(1/3) - T_0^p(1/3)\big)$.
From this equality we find that $S_p(1/3)=2^p/(6^p-3^p)$.

3) Substituting the values $x=1/5$ and $x=2/5$ into the formula~\eqref{eqFuncUrM1}, we get the following two equalities:
the first is $S_p(2/5) = 2^p\big(S_p(1/5) - T_0^p(1/5)\big)$ 
and the second is $S_p(4/5) =2^p\big(S_p(2/5) - T_0^p(2/5)\big)$.

By virtue of~\eqref{eqDefT0} we have: $T_0^p(1/5)=1/5$ and $T_0^p(2/5)=2/5$.
In addition, we have $S_p(4/5) = S_p(1/5)$ due to the symmetry property~\eqref{eqSimmSp}.
With this in mind, from the previous two equalities we get the following two equalities:
$S_p(2/5) = 2^p\big(S_p(1/5) -(1/5)^p\big)$ and
$S_p(1/5) = 2^p\big(S_p(2/5) -(2/5)^p\big)$.
Solving the resulting system of two equations with respect to $S_p(1/5)$ and $S_p(2/5)$, we find the result:
$S_p(1/5)=2\cdot4^p/(20^p-5^p)$ and $S_p(2/5)=(8^p+2^p)/(20^p-5^p)$.
\end{proof}

\section{Global maxima of Takagi power class functions for parameter values $0<p<1$}
\label{sectMaxSp}

This section is devoted to finding the global maxima of the functions $S_p(x)$ for parameter values $p\in(0;1)$, and the sets where they are reached. Moreover, here we prove that for any $p\in(0;1)$ the function $S_p(x)$ has no derivative at points $x=1/3$ and $x=2/3$.
\begin{lemma} \label{LemSab}
For any real numbers $a,b,s$ and $p$, satisfying the conditions $0<s\leq a<b$ and $0<p<1$, the inequality holds
$$
(a+s)^ p-(a-s)^ p > (b+s)^ p-(b-s)^ p.
$$
\end{lemma}
\begin{proof}
Let's define the function $f(s) = \big((a+s)^ p-(a-s)^p\big) - \big((b+s)^ p-(b-s)^ p\big)$ on the segment $[0,a]$.
It is enough for us to prove the inequality $f(s)>0$ for any $s\in(0,a]$.

Calculate the derivative of $f'(s)$ for $s\in(0,a]$:
$$
f'(s) = p\cdot\big((a+s)^{ p-1}-(b+s)^{p-1} + (a-s)^{ p-1}-(b-s)^{ p-1}\big).
$$
Since $(a+s)^{ p-1}>(b+s)^{p-1}$ and $(a-s)^{p-1}>(b-s)^{p-1}$ for all $s\in(0,a]$,
then $f'(s)>0$.
Since the function $f(s)$ is continuous at $s\in[0,a]$, then $f(s)$ strictly increases on the segment $[0,a]$.
By virtue of the equality $f(0)=0$, this implies the inequality $f(s)>0$ for any $s\in(0,a]$.
Which was to be proved.
\end{proof}

\begin{lemma} \label{LemDn}
Let $0< p<1$. For each non-negative integer $n$, let us define the following function on the half-interval $(0;1]$ 
\begin{equation} \label{eqDn}
D_n(s) = \sum_{k=0}^n \Big( \big(2^{n-k+2}-(-1)^{n+k}+(-1)^k\cdot 3s\big)^ p - \big(2^{n-k+2}-(-1)^{n+k}-(-1)^k\cdot 3s\big)^ p \Big).
\end{equation}
Then for any $s\in(0;1]$, the inequality $D_n(s)>0$ will be fulfilled if $n$ is even, and the inequality $D_n(s)<0$ if $n$ is odd.
\end{lemma}
\begin{proof}
For each integer $k=0,1,\ldots,n$ and any $s\in(0;1]$ we put
\begin{equation} \label{eqdnk}
d_{n,k}(s) = \big(2^{n-k+2}-(-1)^{n+k}+(-1)^k\cdot 3s\big)^ p - \big(2^{n-k+2}-(-1)^{n+k}-(-1)^k\cdot 3s\big)^ p.
\end{equation}
Then from the formula~\eqref{eqDn} follows the equality $D_n(s) = \sum_{k=0}^n d_{n,k}(s)$.

1) First we prove that $D_n(s)<0$ for $n$ odd.
To do this, divide the sum $\sum_{k=0}^n d_{n,k}(s)$ defining $D_n(s)$, on pairs of terms:
$D_n(s) = \sum_{i=0}^{(n-1)/2} \big(d_{n,2i}(s)+d_{n,2i+1}(s)\big)$.
It is enough for us to show that the sum of each pair is negative, that is, $d_{n,2i}(s)+d_{n,2i+1}(s)<0$ for any $s\in(0;1]$ and any $i=0,1,\ldots,(n-1)/2$.

Let's put $a_{n,i}=2^{n-2i+1}-1$ and $b_{n,i}=2^{n-2i+2}+1$. 
Then, by virtue of the formula~\eqref{eqdnk}, the equality is true
$$
d_{n,2i}(s)+d_{n,2i+1}(s) = 
(2^{n-2i+2}+1+3s)^ p - (2^{n-2i+2}+1-3s)^ p + (2^{n-2i+1}-1-3s)^ p -
$$
$$
- (2^{n-2i+1}-1+3s)^ p =
\big( (b_{n,i}+3s)^ p-(b_{n,i}-3s)^ p \big) - \big( (a_{n,i}+3s)^ p-(a_{n,i}-3s)^ p \big).
$$
Since $0<3s\leq3\leq a_{n,i}<b_{n,i}$ for any $s\in(0;1]$ and any $i=0,1,\ldots,(n-1)/2$, then the lemma~\ref{LemSab} can be applied. According to this lemma, the estimate is correct
$(a_{n,i}+3s)^ p-(a_{n,i}-3s)^ p > (b_{n,i}+3s)^ p-(b_{n,i}-3s)^ p$.
From here we get the necessary inequality $d_{n,2i}(s)+d_{n,2i+1}(s)<0$ for any $s\in(0;1]$ and any $i=0,1,\ldots,(n-1)/2$.

2) Now we prove that $D_n(s)>0$ for $n$ even.
To do this, in the sum $\sum_{k=0}^n d_{n,k}(s)$, defining $D_n(s)$, we write the first term separately, and divide the remaining terms into pairs:
$D_n(s) = d_{n,0}(s) + \sum_{i=1}^{n/2} \big(d_{n,2i-1}(s)+d_{n,2i}(s)\big)$.
It is enough for us to show that here all terms are positive, that is, $d_{n,0}(s)>0$ and $d_{n,2i-1}(s)+d_{n,2i}(s)>0$ for any $s\in(0;1]$ and any $i=0,1,\ldots,n/2$.

By virtue of the equality~\eqref{eqdnk} we have: $d_{n,0}(s) = (2^{n+2}-1+3s)^ p - (2^{n+2}-1-3s)^ p > 0$.

Next, let's put $a_{n,i}=2^{n-2i+2}-1$ and $b_{n,i}=2^{n-2i+3}+1$.
Then the equality is true
$$
d_{n,2i-1}(s)+d_{n,2i}(s) =
(2^{n-2i+3}+1-3s)^ p - (2^{n-2i+3}+1+3s)^ p + (2^{n-2i+2}-1+3s)^ p -
$$
$$
- (2^{n-2i+2}-1-3s)^ p =
\big( (a_{n,i}+3s)^ p-(a_{n,i}-3s)^ p \big) - \big( (b_{n,i}+3s)^ p-(b_{n,i}-3s)^ p \big).
$$
Since $0<3s\leq3\leq a_{n,i}<b_{n,i}$ for any $s\in(0;1]$ and any $i=0,1,\ldots,n/2$, then the lemma~\ref{LemSab} can be applied. According to this lemma, the following estimate is correct:
$(a_{n,i}+3s)^ p-(a_{n,i}-3s)^ p > (b_{n,i}+3s)^ p-(b_{n,i}-3s)^ p$.
From here we get the necessary inequality $d_{n,2i-1}(s)+d_{n,2i}(s)>0$ for any $s\in(0;1]$ and any $i=0,1,\ldots,n/2$.

The lemma is proved.
\end{proof}

\begin{lemma} \label{lemAnBn}
For each $n=0,1,2,\ldots$ let us define a segment $[a_n,b_n]$ with endpoints
\begin{equation} \label{eqAnBn}
a_n = \frac13+\frac{(-1)^{n+1}}{3\cdot2^{n+2}}-\frac1{2^{n+2}} \quad\text{and}\quad
b_n = \frac13+\frac{(-1)^{n+1}}{3\cdot2^{n+2}}+\frac1{2^{n+2}}.
\end{equation}
Then

1) For each $n=0,1,2,\ldots$ the length of the segment $[a_n,b_n]$ equals to $1/2^{n+1}$, and its middle is the point $c_n = 1/3+(-1)^{n+1}/(3\cdot2^{n+2})$, that can also be written as $c_n=q_n/2^{n+2}$, where $q_n$ is some integer.

2) If $n$ is even, then $a_{n+1}=c_n$ and $b_{n+1}=b_n$, that is, the segment $[a_{n+1},b_{n+1}]$ is the right half of the segment $[a_n,b_n]$.
If $n$ is odd, then $a_{n+1}=a_n$ and $b_{n+1}=c_n$, that is, the segment $[a_{n+1},b_{n+1}]$ is the left half of the segment $[a_n,b_n]$.

3) The segments $[a_n,b_n]$ are nested (that is
$[0,1/2]=[a_0,b_0]\supset[a_1,b_1]\supset\ldots\supset[a_n,b_n]\supset[a_{n+1},b_{n+1}]\supset\ldots$), 
and the equality $\bigcap_{n=0}^\infty[a_n,b_n]=\{1/3\}$ holds.
\end{lemma}
\begin{proof}
1) At this item, the only thing that needs proof is that the number $c_n = 1/3 +(-1)^{n+1}/(3\cdot2^{n+2})$ can be represented as $c_n=q_n/2^{n+2}$, where $q_n$ is an integer. 
From the equality $c_n = 1/3+(-1)^{n+1}/(3\cdot2^{n+2}) = q_n/2^{n+2}$ we find that $q_n = (2^{n+2}+(-1)^{n+1})/3$.
If $n$ is even, then 
the remainder in division $2^{n+2}$ by $3$ equals to $1$, and $(-1)^{n+1}=-1$. So $q_n$ is an integer.
Similarly, if $n$ is odd, then 
the remainder in division $2^{n+2}$ by $3$ equals to $2$, and $(-1)^{n+1}=1$. So in this case $q_n$ is an integer too.

2) If $n$ is even, then the following four equalities follow from the formula~\eqref{eqAnBn}:

$c_n = 1/3-1/(3\cdot2^{n+2})$, 

$b_n = 1/3-1/(3\cdot2^{n+2})+1/2^{n+2} = 1/3+1/(3\cdot2^{n+1})$,

$a_{n+1} = 1/3+1/(3\cdot2^{n+3})-1/2^{n+3} = 1/3-1/(3\cdot2^{n+2})$ and

$b_{n+1} = 1/3+1/(3\cdot2^{n+3})+1/2^{n+3} = 1/3+1/(3\cdot2^{n+1})$.

Thus, for $n$ even we have: $a_{n+1}=c_n$ and $b_{n+1}=b_n$.
For $n$ odd, check of the equality $a_{n+1}=a_n$ and $b_{n+1}=c_n$ is done similarly.

3) The nesting of the segments $[a_n,b_n]$ follows from item 2 of the statement of this lemma.
The equality $\bigcap_{n=0}^\infty[a_n,b_n]=\{1/3\}$ follows from Cantor's theorem on~nested segments
and from the relation 

$\displaystyle\lim_{n\to\infty}a_n = \lim_{n\to\infty}\big(1/3+(-1)^{n+1}/(3\cdot2^{n+2})-1/2^{n+2}\big) = 1/3$.
\end{proof}

\begin{theorem} \label{teorMaxSp}
For any $p\in(0;1)$, the global maximum of the function $S_p(x)$ on the segment $[0;1]$ is reached only at the points $x=1/3$ and $x=2/3$, and its value is equal to $\max_{x\in[0;1]}S_p(x)=S_p(1/3)=2^p/(6^p-3^p)$.
\end{theorem}
\begin{proof}
1) First we prove that the set of points of the global maximum of the function $S_p(x)$ on the segment $[0;1]$ consists only of points $x=1/3$ and $x=2/3$.

Let's denote the set of points of the global maximum of the function $S_p(x)$ on the segment $[0;1/2]$ by $Argmax_{[0,1/2]}S_p$.
By virtue of the symmetry property~\eqref{eqSimmSp} we have: $S_p(x)=S_p(1-x)$ for all $x\in[0;1]$.
Therefore, it is enough for us to prove that $Argmax_{[0,1/2]}S_p = \{1/3\}$.
To achieve this goal, we show that the set $Argmax_{[0,1/2]}S_p$ is contained in each of the segments $[a_n,b_n]$ 
($n=0,1,2,\ldots$) given in lemma~\ref{lemAnBn}.
Hence, by virtue of point~3 of Lemma~\ref{lemAnBn}, the relation $Argmax_{[0,1/2]}S_p\subset\bigcap_{n=0}^\infty[a_n,b_n]=\{1/3\}$ will follow, from which the desired equality $Argmax_{[0,1/2]}S_p = \{1/3\}$ will be obtained.

The inclusion $Argmax_{[0,1/2]}S_p\subset[a_n,b_n]$ for each $n=0,1,2,\ldots$ can be proved by induction.

\emph{Induction base}: $n=0$.
The inclusion $Argmax_{[0,1/2]}S_p\subset[a_0,b_0]=[0,1/2]$ is true by definition of the set $Argmax_{[0,1/2]}S_p$.

\emph{Induction step}: provided that the inclusion $Argmax_{[0,1/2]}S_p\subset[a_n,b_n]$ is true, it is necessary to prove the inclusion $Argmax_{[0,1/2]}S_p\subset[a_{n+1},b_{n+1}]$.
To do this, we use the following function on the half-interval $(0;1]$: 
\begin{equation} \label{eqOprFn}
f_n(s) = S_p(c_n+s/2^{n+2})-S_p(c_n-s/2^{n+2}),
\end{equation}
where $c_n$ is the middle of the segment $[a_n,b_n]$ (see lemma~\ref{lemAnBn}).
Note that if the point $s$ runs through the half-interval $(0;1]$, then the point $c_n+s/2^{n+2}$ runs through the right half of the segment $[a_n,b_n]$ (except its middle), and the point $c_n-s/2^{n+2}$ runs through the left half of the segment $[a_n,b_n]$ (also except its middle).
Therefore, if for all $s\in(0;1]$ the inequality $f_n(s)>0$ is true, then the set $Argmax_{[0,1/2]}S_p$ lies in the right half of the segment $[a_n,b_n]$, and if for all $s\in(0;1]$ the inequality $f_n(s)<0$ is true, then the set $Argmax_{[0,1/2]}S_p$ lies in the left half of the segment $[a_n,b_n] $ (from the further line of argument it will be seen that only these two options are possible).
Therefore, it is necessary to investigate the sign of the function $f_n(s)$ at $s\in(0;1]$.
By virtue of the formulas \eqref{eqOprFn} and~\eqref{eqDefSp} we have:
\begin{equation} \label{eqFnPreobr}
\begin{split}
f_n(s) = \sum_{k=0}^\infty 
\bigg( \frac{T_0\big(2^k(c_n+s/2^{n+2})\big)}{2^k}\bigg)^p - 
\bigg( \frac{T_0\big(2^k(c_n+s/2^{n+2})\big)}{2^k}\bigg)^p =\\
= \sum_{k=0}^\infty \frac1{2^{kp}}
\Big(T_0^p(2^k c_n + 2^{k-n-2}s) - T_0^p(2^k c_n + 2^{k-n-2}s)\Big). 
\end{split}
\end{equation}
Let's use the symmetry of the function $T_0$ with respect to half-integers, that is, numbers of the form $q/2$, where $q$ is any integer number (see~\eqref{eqSimmSp}).
By virtue of item~1 of Lemma~\ref{lemAnBn}, the number $c_n$ can be written as $c_n=q_n/2^{n+2}$, where $q_n$ is an integer.
Therefore, $2^k c_n=2^kq_n/2^{n+2}=2^{k-n-1}q_n/2$, and for $k\geq n+1$, the number $2^k c_n$ will be a half-integer.
Therefore, for $k\geq n+1$, the function $T_0$ will be symmetric with respect to the points $2^k c_n$, that is, the equality $T_0(2^k c_n +2^{k-n-2}s)=T_0(2^k c_n- 2^{k-n-2}s)$ will be fulfilled.
Therefore, on the right side of the equality~\eqref{eqFnPreobr}, all terms with numbers $k\geq n+1$ are equal to zero.
Thus, the formula~\eqref{eqFnPreobr}, by virtue of the equality $c_n = 1/3+(-1)^{n+1}/(3\cdot2^{n+2})$ from item~1 of lemma~\ref{lemAnBn}, can be rewritten as
\begin{equation} \label{eqFnPreobr2}
\begin{split}
f_n(s) = \sum_{k=0}^n \frac1{2^{kp}}
\Big(T_0^p(2^k c_n + 2^{k-n-2}s) - T_0^p(2^k c_n + 2^{k-n-2}s)\Big) =\\
= \sum_{k=0}^n \frac1{2^{kp}}\bigg(
T_0^p\Big(\frac{2^k}{3} + \frac{(-1)^{n+1}+3s}{3\cdot 2^{n-k+2}}\Big) - 
T_0^p\Big(\frac{2^k}{3} + \frac{(-1)^{n+1}-3s}{3\cdot 2^{n-k+2}}\Big) \bigg).
\end{split}
\end{equation}
Due to the function $T_0$ has the period $1$, in this equality, numbers of the form $2^k/3$ can be replaced by their fractional parts $\{2^k/3\}$.
If $k\geq0$ and is even, then $\{2^k/3\}=1/3$, so $\{2^k/3\}+((-1)^{n+1}\pm3s)/(3\cdot 2^{n-k+2})\in[0,1/2]$ and, by virtue of the formula~\eqref{eqDefT0} for the function $T_0$, the equalities are true
\begin{equation} \label{eqT0kCet}
T_0\Big(\frac{2^k}{3} + \frac{(-1)^{n+1} \pm 3s}{3\cdot 2^{n-k+2}}\Big) = 
\frac13 + \frac{(-1)^{n+1}\pm 3s}{3\cdot 2^{n-k+2}}\quad\text{for even }k\geq0.
\end{equation}
Similarly, if $k\geq0$ and is odd, then $\{2^k/3\}=2/3$, so $\{2^k/3\}+((-1)^{n+1}\pm3s)/(3\cdot 2^{n-k+2})\in[1/2,1]$ and, by virtue of the formula~\eqref{eqDefT0}, the equalities are true
\begin{equation} \label{eqT0kNec}
\begin{split}
T_0\Big(\frac{2^k}{3} + \frac{(-1)^{n+1} \pm 3s}{3\cdot 2^{n-k+2}}\Big) &= 
T_0\Big(\Big\{\frac{2^k}{3}\Big\} + \frac{(-1)^{n+1} \pm 3s}{3\cdot 2^{n-k+2}}\Big) =\\
= 1- \Big(\frac23 + \frac{(-1)^{n+1} \pm 3s}{3\cdot 2^{n-k+2}}\Big) &= 
\frac13 - \frac{(-1)^{n+1}\pm 3s}{3\cdot 2^{n-k+2}}\quad\text{for odd}k\geq0.
\end{split}
\end{equation}
Given the formulas \eqref{eqT0kCet} and~\eqref{eqT0kNec}, we conclude that for all integers $k\geq0$ the following equality is true:
\begin{equation*} 
T_0\Big(\frac{2^k}{3} + \frac{(-1)^{n+1} \pm 3s}{3\cdot 2^{n-k+2}}\Big) = 
\frac13 + (-1)^k\frac{(-1)^{n+1} \pm 3s}{3\cdot 2^{n-k+2}}.
\end{equation*}
Hence, the following chain of equalities follows from the formula~\eqref{eqFnPreobr2}:
\begin{equation*} 
\begin{split}
f_n(s) &= \sum_{k=0}^n \frac1{2^{kp}}\bigg(
\Big(\frac13 + (-1)^k\frac{(-1)^{n+1} + 3s}{3\cdot 2^{n-k+2}}\Big)^p - 
\Big(\frac13 + (-1)^k\frac{(-1)^{n+1} - 3s}{3\cdot 2^{n-k+2}}\Big)^p\bigg) =\\
&= \frac{1}{3^p\cdot2^{(n+2)p}} \sum_{k=0}^n
\Big( \big(2^{n-k+2}-(-1)^{n+k}+(-1)^k\cdot 3s\big)^p -\\
&\qquad\qquad\qquad\qquad- \big(2^{n-k+2}-(-1)^{n+k}-(-1)^k\cdot 3s\big)^p \Big) =
\frac{D_n(s)}{3^p\cdot2^{(n+2)p}},
\end{split}
\end{equation*}
where $D_n(s)$ is a function defined by the equality~\eqref{eqDn} and studied in the lemma~\ref{LemDn}.

If $n$ is even, then, by virtue of the lemma~\ref{LemDn}, the function $D_n(s)$, and hence the function $f_n(s)$, is positive for $s\in(0;1]$, so the set $Argmax_{[0,1/2]}S_p$ lies in the right half of the segment $[a_n,b_n]$, that is, on the segment $[a_{n+1},b_{n+1}]$ (according to point~4 lemmas~\ref{lemAnBn}).
Similarly, if $n$ is odd, then, again by virtue of the lemma~\ref{LemDn}, the function $D_n(s)$, and hence the function $f_n(s)$, is negative for $s\in(0;1]$, so the set $Argmax_{[0,1/2]}S_p$ lies in the left half of the segment $[a_n,b_n]$, that is, again on the segment $[a_{n+1},b_{n+1} ]$ (according to item~4 of the lemma~\ref{lemAnBn}).
Thus, the inclusion of $Argmax_{[0,1/2]}S_p\subset[a_{n+1},b_{n+1}]$ takes place for both cases: $n$ is even or $n$ is odd.

So, the induction step is done, and the first statement of the theorem is proved.

2) Let us calculate the value of the global maximum of the function $S_p(x)$ on the segment $[0;1]$.
By virtue of the first item of this proof, this value is equal to the value of the function $S_p$ at the point $x=1/3$, which, by virtue of the lemma~\ref{lemSp121315}, is equal to $2^p/(6^p-3^p)$.

This completes the proof.
\end{proof}

\begin{theorem}
For any real $p\in(0;1)$, the function $S_p(x)$ has no derivative at the points $x=1/3$ and $x=2/3$.
\end{theorem}
\begin{proof}
1) At first we prove by contradiction that $S_p(x)$ is not differentiable at~point $x=1/3$.
Suppose that the derivative $S_p'(1/3)$ exists.
Then, differentiating the functional equation~\eqref{eqFuncUrM1} at the point $x=1/3$, we obtain that the derivative $S_p'(2/3)$ also exists, and the equality holds
$2S_p'(2/3) = 2^p\big(S_p'(1/3) - pT_0^{p-1}(1/3)T_0'(1/3)\big)$.
Hence, considering that $T_0(1/3)=1/3$ and $T_0'(1/3)=1$, we get the formula
\begin{equation} \label{eqSpDif13}
2S_p'(2/3) = 2^p\big(S_p'(1/3) - p\cdot 3^{1-p}\big).
\end{equation}
Since, by virtue of the theorem~\ref{teorMaxSp}, the points $x=1/3$ and $x=2/3$ are extremum points, the equalities $S_p'(1/3) = S_p'(2/3) = 0$ must be satisfied.
Substituting these values into the formula~\eqref{eqSpDif13}, we get an incorrect equality $0 = -2^p p\cdot 3^{1-p}$.

The resulting contradiction shows that $S_p(x)$ is not differentiable at~point $x=1/3$.

2) Since the derivative $S_p'(1/3)$ does not exist, then, by virtue of the equality~\eqref{eqSpDif13}, the derivative $S_p'(2/3)$ does not exist either.
\end{proof}
\section{Acknowledgements}

The publication was prepared within the framework of the Academic Fund Program at HSE University in 2020–2021 (grant No.20-04-022, project ``Evolution semigroups and their applications'').
The authors thank the head of this project, Ivan Remizov, for his attention to the work.

\end{document}